\documentclass[11pt,a4paper,reqno]{amsart}
\usepackage[english]{babel}
\usepackage[applemac]{inputenc}
\usepackage{bbm}
\usepackage{cite}
\usepackage{float, verbatim}
\usepackage{amsmath}
\usepackage{amssymb}
\usepackage{amsthm}
\usepackage{amsfonts}
\usepackage{graphicx}
\usepackage{mathtools}
\usepackage{enumitem}
\usepackage[colorlinks = true, citecolor = black]{hyperref}
\usepackage{color}
\usepackage{tikz}
\usepackage{xcolor}
\usepackage{pgfplots}
\usepackage{caption}
\usepackage{subcaption}
\usepackage{enumerate}
\usepackage{autobreak}
\usepackage{url}

\newcommand{\eps}{\varepsilon}
\renewcommand{\epsilon}{\varepsilon}

\newtheorem*{thm*}{Theorem}
\newtheorem*{conj*}{Conjecture}

\newtheorem{thm}{Theorem}[section]

\newtheorem{conj}[thm]{Conjecture}

\newtheorem{ques}[thm]{Question}

\begin{document}
	
	\title{Approximate arithmetic structure in large sets of integers}
	
	\author{Jonathan M. Fraser}
	\address{J. M. Fraser\\
		School of Mathematics \& Statistics\\University of St Andrews\\ St Andrews\\ KY16 9SS\\ UK  }
	\curraddr{}
	\email{jmf32@st-andrews.ac.uk}
	\thanks{}

\author{Han Yu}
	\address{H. Yu\\
		School of Mathematics \& Statistics\\University of St Andrews\\ St Andrews\\ KY16 9SS\\ UK  }
	\curraddr{}
	\email{hy25@st-andrews.ac.uk}
	\thanks{}
	
	
	
	\date{}
	
	\dedicatory{}
	
	\begin{abstract}
	We prove that if a set is `large' in the sense of Erd\H{o}s, then it approximates arbitrarily long arithmetic progressions in a strong quantitative sense. More specifically, expressing the error in the approximation in terms of the gap length $\Delta$ of the progression, we improve a previous result of $o(\Delta)$ to $O(\Delta^\alpha)$ for any $\alpha \in (0,1)$.
	\end{abstract}
	
	\maketitle
	\allowdisplaybreaks

   \section{The Erd\H{o}s conjecture on arithmetic progressions}

The  Erd\H{o}s conjecture on arithmetic progressions is a famous open problem in combinatorial number theory, which states that a set of positive integers whose reciprocals form a divergent series\footnote{ Erd\H{o}s used the term `large' to describe such sets.  } should contain arbitrarily long arithmetic progressions.  More precisely, if    $A\subset\mathbb{N}$ is such that  $\sum_{a\in A} a^{-1}=\infty$, then $A$ should contain an arithmetic progression of length $k$ for all $k \geq 1$. Note that  the Green-Tao theorem, where $A$ is the set of prime numbers, is a special case of this conjecture, see \cite{greentao}.  

In \cite{fraseryu}, we established the following weak version of this conjecture.  We say an arithmetic progression $\{a, a+\Delta, \dots, a +(k-1)\Delta\} \subset \mathbb{R}$ is of \emph{length} $k \in \mathbb{N}$ and \emph{gap length} $\Delta>0$.

\begin{thm}[Theorem 2.1 in \cite{fraseryu}] \label{ap1}
If    $A\subset\mathbb{N}$ is such that  $\sum_{a\in A} a^{-1}=\infty$, then for all $\eps >0$ and all $k \geq 1$,   there exists an arithmetic progression $P$ of length $k$ and gap length $\Delta \geq 1$ such that
\[
\sup_{p \in P} \inf_{a \in A} |p-a| \leq \eps \Delta.
\]
\end{thm}

This result should be interpreted as saying that $A$ gets \emph{arbitrarily close} to arbitrarily long arithmetic progressions, see also \cite{fraserAP}. One disadvantage of this result is that the conclusion holds for natural examples of sets which are known \emph{not} to contain long arithmetic progressions, such as the squares.  Moreover, there is no information on how $P$ depends on $\eps$ and $k$. Another direction in which this result could be improved is to allow $\eps$ to depend on $\Delta$, that is to improve the bound on the `uncertainty' from $o(\Delta)$ to something stronger.    For example, if one could show that $\epsilon$ could be replaced by $C\Delta^{-1}$ for some constant $C,$ then the Erd\H{o}s  conjecture would follow. This is because our `uncertainty' $\eps \Delta$ is at most $C$ and the result would then follow from  Van der Waerden's Theorem, see \cite{VW}.  More precisely, let $P$ be a very long  arithmetic progression which our set approximates to within $C$.  Colour the points in our set which approximate $P$ according to their distance from $P$ (in particular we need at most 2$C$+1 colours).  Van der Waerden's Theorem allows us to extract a monochromatic arithmetic progression of length $k$ provided the length of $P$ is large enough.  Since this progression is monochromatic it is a genuine arithmetic progression of length $k$ inside our set.

The main result of this paper addresses each of the above issues.  For example, the conclusion of our main result should not hold for sets such as the squares or cubes (although proving this rigorously seems challenging and is related to Mazur's \emph{near miss problem}, see Section \ref{Mazur}), and we get quantitative  information regarding  $P$.  We essentially show  that one can choose $\eps=C\Delta^{-\delta}$ for any $0<\delta<1.$ Thus, we push Theorem \ref{ap1} further towards  Erd\H{o}s  conjecture, which would follow if we could choose  $\delta=1$.    The ideas in this paper can be used to push the result even further, relying on a result of Gowers \cite{Go}, see Section \ref{Fur}. 

\section{Main result}

Our main result is the following improvement over Theorem \ref{ap1}. We write $\# E$ to denote the cardinality of a finite set $E$.

   \begin{thm}\label{Main}
Suppose $A\subset\mathbb{N}$ is such that  there exists a constant $\gamma>0$ such that
\[
\#A \cap [0,n]  \geq \frac{n}{(\log n)^\gamma}
\]
for infinitely many $n$.  Then,  for all $\alpha\in (0,1), k \geq 1, \Delta_0>1$, there exists  infinitely many arithmetic progressions $P$ of length $k$ and gap length $\Delta \geq \Delta_0$ such that
\[
\sup_{p \in P} \inf_{a \in A} |p-a| \leq \Delta^\alpha.
\]
Moreover, there is a constant $c>0$, depending only on $\alpha, \gamma$, such that for infinitely many $n \in\mathbb{N}$,  $P$  can be chosen to have gap length at least $cn$ and lie in the interval $\big[2^n, \, 2^{n+1}\big]$.
   \end{thm}

Since sets of integers whose reciprocals form a divergent series necessarily  satisfy the power-log density assumption above (with  $\gamma>1$), Theorem \ref{Main} applies to sets which are `large' in the sense of Erd\H{o}s.  More precisely, if   $A\subset\mathbb{N}$ is such that  $\sum_{a\in A} a^{-1}=\infty$, then for all $\gamma>1$, we have
\[
\#A \cap [0,n]  \geq \frac{n}{(\log n)^\gamma}\tag{*}
\]
for infinitely many $n$.  Indeed, if (*) were not true for some $\gamma>1$ and all $n \geq e^{t}$ for some positive integer $t$, then
\[
\sum_{a\in A} a^{-1} \  \leq \ \sum_{a \in A \cap [1, e^t]} \frac{1}{a} \,  + \,   \sum_{n \geq t} \sum_{a \in A \cap [e^n,e^{n+1}]} e^{-n}\   \leq \  2 t \, + \, e\sum_{n \geq 1}  n^{-\gamma} \ < \ \infty.
\]
 Finally, we note that Theorem \ref{ap1} follows directly from Theorem \ref{Main} as follows.  Fix $\eps >0$ and  $k \geq 1$ and apply  Theorem \ref{Main} with $\alpha=1/2$ and $\Delta$ chosen  large enough to ensure that $\sqrt{\Delta} \leq \eps \Delta$.

\section{Proof of Theorem \ref{Main}}\label{PROOF}

Fix $\alpha\in (0,1)$, $k \geq 3$ and $\gamma>0$. Let $\epsilon \in (0,1/2)$   and $M_\epsilon$ be an integer such that, for all $M>M_\eps$, any subset of $\{1,\dots,M\}$ with at least $\epsilon M$ elements must contain an arithmetic progression of length $k$.  Such a number $M_\epsilon$ exists as a direct consequence of Szemer\'{e}di's celebrated theorem on arithmetic progressions, \cite{SZ}.    Later we will fix a particular $\eps$ depending on $\alpha$ and $\gamma$.

  We write $X \lesssim Y$ to mean that $X \leq cY$ for some universal constant $c>0$.  We also write $X \gtrsim Y$ to mean $Y \lesssim X$ and $X \approx Y$ to mean that both $X \lesssim Y$ and $X \gtrsim Y$ hold.   The implicit constants $c$ appearing here can, and often will, depend on $\alpha, k, \gamma, \eps$.  We write $\lfloor x \rfloor$ to denote the integer part of $x \geq 0$. 

 Consider the dyadic intervals $[2^n,2^{n+1})$ for integers $n\geq 0$ and let $A_n=A\cap [2^n,2^{n+1}).$ For notational  simplicity we write $N=2^n.$  Decompose $[2^n,2^{n+1})$ into smaller (half-open) intervals of equal length $ N^{\alpha^2}$ and label these intervals from left to right by $1, 2, \dots$.  We may not be able to perform this decomposition exactly, in which case we will be left with an interval of length $N-\lfloor N/N^{\alpha^2} \rfloor N^{\alpha^2}$ at the right hand side, which we simply discard. Group these  intervals into equivalence classes by considering their labels modulo $\lfloor N^\alpha/N^{\alpha^2} \rfloor$.  Note that the set of centres of intervals in a given equivalence class form an arithmetic progression of length (at least) $ \tfrac{1}{2} N/N^\alpha  $ and gap length  $N^{\alpha^2} \lfloor N^\alpha/N^{\alpha^2} \rfloor$.  

Suppose that $\tfrac{1}{2} N/N^{\alpha}>M_\epsilon$, which is certainly true for all sufficiently large $n$.  In  order to avoid the existence of an arithmetic progression $P$  with length $k$ and gap  length $\Delta =  N^{\alpha^2} \lfloor N^\alpha/N^{\alpha^2} \rfloor$ such that
\begin{equation} \label{222}
\sup_{p \in P} \inf_{a \in A_n} |p-a| \leq 2 \Delta^\alpha,
\end{equation}
we see that $A_n$ can intersect no more than $  2\epsilon  N/N^{\alpha^2}  $ many intervals of length $N^{\alpha^2}.$   Indeed, each equivalence class contains at least $M_\eps$ and at most $2 N/N^\alpha $ many intervals and so must intersect $A_n$ in fewer than $ 2 \eps  N/N^\alpha $ many of these intervals, and there are fewer than $N^\alpha/N^{\alpha^2} $ many equivalence classes. 

We repeat the above  argument inside each  interval of length $ N^{\alpha^2}$ that intersects $A_n.$  That is, we decompose each interval of length $N^{\alpha^2}$ into intervals of equal length $N^{\alpha^4}$ (possibly discarding a small piece at the end) and then we work modulo modulo $\lfloor N^{\alpha^3}/N^{\alpha^4} \rfloor$.   In particular, if $\tfrac{1}{2}  N^{\alpha^2}/N^{\alpha^3}  >M_\epsilon$ and there does not exist an arithmetic progression $P$  with length $k$ and gap  length $\Delta =  N^{\alpha^4} \lfloor N^{\alpha^3}/N^{\alpha^4} \rfloor$ satisfying \eqref{222}, then $A_n$ can intersect no more than  
   \[
  (2\epsilon)^2 \frac{N}{N^{\alpha^2}}\frac{N^{\alpha^2}}{N^{\alpha^4}}
   \]
   many intervals of length $N^{\alpha^4}.$ We can repeat this decomposition argument $(m'+1)$ times where $m'$ is chosen to be the largest integer satisfying
\[
\tfrac{1}{2} N^{\alpha^{2m'}}/N^{\alpha^{2m'+1}}   >M_\epsilon,
\]
noting that
   \[
   m' \approx \frac{\log\log N + \log (1-\alpha)-\log\log M_\epsilon}{-2\log \alpha}.\tag{**}
   \]
However, repeating the argument this many times means we can only bound the gap length of the arithmetic progressions we avoid below by
\[
\Delta \geq N^{\alpha^{2m'+2}} \lfloor N^{\alpha^{2m'+1}}/N^{\alpha^{2m'+2}} \rfloor \geq \tfrac{1}{2}N^{\alpha^{2m'+1}} \approx 1
\]
which is not sufficient to prove the theorem.  Therefore we choose to repeat the argument only $(m+1)$ times where
\[
m=m'-\frac{\log\log\log N}{-2 \log \alpha}.
\]
At the $(l+1)$st step we decompose intervals of length $N^{\alpha^{2l}}$ into intervals of length $N^{\alpha^{2l+2}}$ and work modulo  $\lfloor N^{\alpha^{2l+1}}/N^{\alpha^{2l+2}} \rfloor$. 

After applying this argument $(m+1)$ times we see that if there  does not exist an arithmetic progression  $P$ satisfying \eqref{222} with gap length 
\[
\Delta \geq \tfrac{1}{2}N^{\alpha^{2m+1}} \gtrsim n,
\]
then $A_n$ contains at most
\[
    (2\epsilon)^{m+1} \frac{N}{N^{\alpha^2}}\frac{N^{\alpha^2}}{N^{\alpha^4}}\cdots \frac{N^{\alpha^{2m}}}{N^{\alpha^{2m+2}}} \left( 2 N^{\alpha^{2m+2}} \right) \ = \ 2(2\epsilon)^{m+1} N 
\]
many elements, where we used the fact that intervals of length $N^{\alpha^{2m+2}}$ contain at most $N^{\alpha^{2m+2}}+1 \leq 2 N^{\alpha^{2m+2}}$ many integers.  In particular, 
\[
\# A_n \leq 2(2\epsilon)^{m+1} N \lesssim  (2\eps)^m 2^n  \lesssim       \left(\frac{\log n}{n}\right)^{\frac{\log 2\eps}{2 \log \alpha}}2^n 
\]
where the implicit constants here depend on   $\alpha$ and $\eps$.  

In order to reach a contradiction,  suppose that for all but finitely many $n$ there does \emph{not} exist an arithmetic progression  $P$ satisfying \eqref{222}. Therefore, the  above cardinality estimate for $A_n$ holds for all but finitely many $n$.  We now fix   $\eps>0$ depending on $\alpha$ and $\gamma$  such that
\[
\left(\frac{\log n}{n}\right)^{\frac{\log 2\eps}{2 \log \alpha}} \leq n^{-\gamma}
\]
for all $n \geq 1$. Therefore,   for integer $T>0$,
\[
\# A \cap [0,T] \leq \sum_{n=0}^{\lceil \frac{\log T}{\log 2}\rceil}  \# A_n  \lesssim    \sum_{n=0}^{\lceil \frac{\log T}{\log 2}\rceil}    n^{-\gamma} 2^n \lesssim    \frac{T}{(\log T)^\gamma}.
\]
This contradicts the power-log density hypothesis since $\gamma>0$ can be chosen to be arbitrarily large.  Therefore, for infinitely many $n$, there exists an arithmetic progression  $P \subset [2^n,2^{n+1})$  satisfying \eqref{222} with gap length $\Delta \gtrsim n$, proving the theorem. Note that  the upper bound of $2 \Delta^\alpha$  in  \eqref{222} can be trivially upgraded to the desired upper bound of $ \Delta^\alpha$ by replacing $\alpha$ with $\alpha' \in (0,\alpha)$ in the argument and choosing $\Delta$ large enough.

\section{A remark on Mazur's near miss problem}\label{Mazur}
From Theorem \ref{Main} we see that all `large enough' sets  $A\subseteq \mathbb{N}$ must `nearly' contain arbitrarily long arithmetic progressions.   It is interesting to consider examples of sets $A$ for which the conclusion of Theorem \ref{ap1} is satisfied but the conclusion of Theorem \ref{Main} fails.  It is easy to show existence of such examples but we are particularly interested in the sets $A_t = \{ n^{-t} : n \in \mathbb{N}\}$ for a fixed integer $t\geq 2$, which turn out to be elusive and related to a deep problem posed by Mazur.  It follows from \cite{fraseryu} that $A_t$ satisfies the conclusion of Theorem \ref{ap1} since the set of reciprocals of elements in $A_t$ is a set of full Assouad dimension.  However, it is known that these sets do not contain genuine arithmetic progressions of length $k \geq 4$ (or $k \geq 3$ provided $t \geq 3$).  We make the following conjecture.  

\begin{conj}\label{Conj}
Let $A_t = \{ n^{-t} : n \in \mathbb{N}\}$ for a fixed integer $t\geq 2$ and let $\alpha \in (0,1)$. There exist  integers $k_0 \geq 3$ and $\Delta_0 \geq 1$ such that if $P$ is an arithmetic progression of length $k \geq k_0$ and gap length  $\Delta \geq \Delta_0$, then 
	\[
	\sup_{p \in P} \inf_{a \in A_t} |p-a| > \Delta^\alpha.
	\] 
\end{conj}

The above conjecture is related to Mazur's \emph{near miss problem}, see \cite[Section 11]{MA}. To illustrate the connection, let us only consider arithmetic  progressions of length 3. If $a,b,c$ forms an arithmetic progression then $a+c=2b.$ Suppose that $a,b,c$ are $t$-th powers, in which case  we can find rational numbers $r$ and $s$ such that
\[
r^t+s^t=2.
\]
Our goal is not finding exact progressions in the set of $t$-th powers - indeed, there are no arithmetic progressions of length 3 inside the cubes. Instead, we want to know how close the set of $t$-th powers can get to an arithmetic progression of length 3. Given a large positive integer $Q$, we are interested in estimating the smallest distance  between points on the lattice $\mathbb{Z}^2/Q$  and the curve  defined by
\[
\{(x,y) \in \mathbb{R}^2 : x^t+y^t=2 \},
\]
which  is in the spirit of Mazur's near miss problem. More specifically, Conjecture \ref{Conj} is related to bounding this distance from below by $Q^{-\alpha}$ for $\alpha\in (0,1).$  From here, it is natural to consider the following question.

\begin{ques}\label{Ques}
Given an integer $t\geq 3$, what is 
\[
f_t : = \inf_{b>a  \geq 10} \frac{\log  \min_{n \in \mathbb{N}} | \frac{a^t+b^t}{2}-n^t | }{\log |b^t-a^t|}?
\]
\begin{comment}
  let $Dio_k(n,m)$ be defined as follows
	\[
	Dio_k(n,m)=|^k((n^k+m^k)/2)-(n^k+m^k)/2|.
	\]
	Then we define $dio_k(n,m)$ to be
	\[
	dio_k(n,m)=\log Dio_k(n,m)/\log |m^k-n^k|.
	\]
	What is the value of
	\[
dio_k=\inf_{n,m\geq 10} dio_k(n,m)? 
	\]
\end{comment}
\end{ques}

\begin{figure}[h]
	\includegraphics[width= 0.9\linewidth, height=10cm]{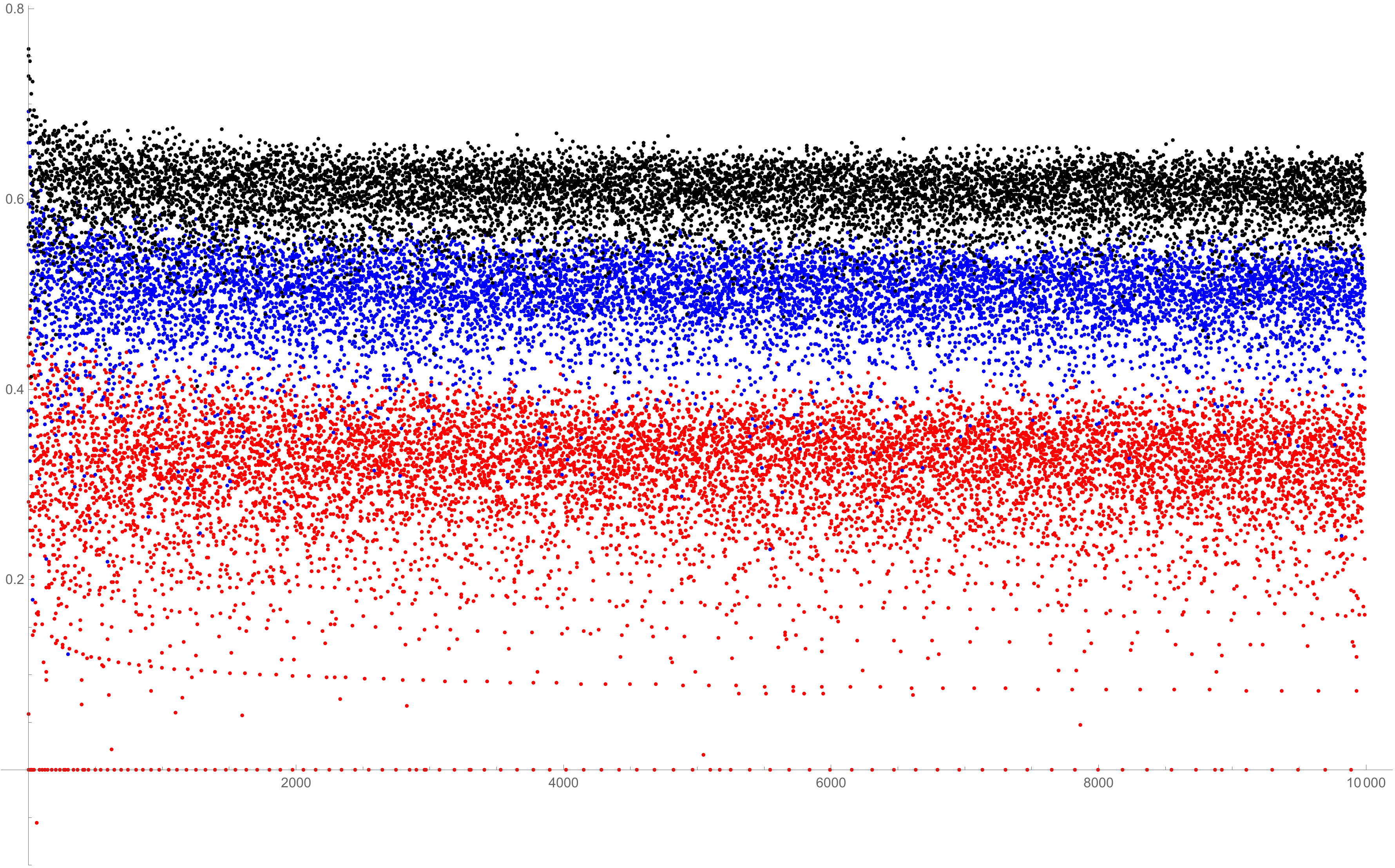} 
	\caption{Plots of $f_t(b)$ as a function of $b$ for different values of $t$.  The plot of $f_3(b)$ is shown in red, $f_4(b)$ is shown in blue, and $f_5(b)$ is shown in black. There is one red point below the $x$-axis which corresponds to $a=42,b=71,n=60.$}
	\label{fig:figure 1}
\end{figure}

In order to test this question numerically, we computed the values
\[
f_t(b) = \min_{b>a  \geq 10} \frac{\log  \min_{n \in \mathbb{N}} | \frac{a^t+b^t}{2}-n^t | }{\log |b^t-a^t|}
\]
for integer values of $b$ up to $10,000$.  The results are plotted in Figure \ref{fig:figure 1}.  This simulation suggests that $f_3 = 0$ but that   $f_t>0$ for $t \geq 4$, and we believe this is the case.
\begin{conj}
	For integers  $t\geq 4$,  $f_t>0.$ On the other hand $f_3=0$ and there are infinitely many integer solutions to
	\[
	x^3+y^3-2z^3 \in \{ \pm 1,\pm 2\}.
	\]
\end{conj}

\section{Further comments}\label{Fur}
\subsection{Higher dimensions}
Szemer\'{e}di's theorem can be generalized for studying `structures' in large subsets of $\mathbb{Z}^d,d\geq 2,$ see \cite{FK78}. This can help us get a higher dimensional version of Theorem \ref{Main} with an almost identical proof as in Section \ref{PROOF}.

%Let $\alpha\in (0,1),\Delta>0$, and $C\subset\mathbb{Z}^d$ be a finite set\footnote{We use $C$ for \textbf{C}onstellation as suggested in \cite{TZ15}.}. 

%We say that $X\subset\mathbb{Z}^{d}$ is an $(\alpha,\Delta)$-approximate of $C$ if there is a translate $v\in\mathbb{Z}^d$ such that the Hausdorff distance between $\Delta C+v$ and  $X$ is at most $\Delta^\alpha.$

\begin{thm}
Let $C\subset\mathbb{Z}^d$ be a finite set\footnote{We use $C$ for \textbf{C}onstellation as suggested in \cite{TZ15}.}. Suppose $A\subset\mathbb{Z}^d$ is such that  there exists a constant $\gamma>0$ such that
	\[
	\#A \cap \{x\in\mathbb{Z}^d:|x|\in [0,n]   \}  \geq \frac{n^d}{(\log n)^\gamma}
	\]
	for infinitely many $n$.  Then,  for all $\alpha\in (0,1), k \geq 1, \Delta_0>1$, there exists  infinitely many  translates $t\in\mathbb{Z}^d$ and $\Delta \geq \Delta_0$ such that 
\[
\sup_{p \in \Delta C+t} \inf_{a \in A} |p-a| \leq \Delta^\alpha,
\]
where $\Delta C+t$ is the image of $C$ under the similitude $z \mapsto \Delta z+t$.
\end{thm}

\subsection{Reducing the uncertainty}
In proving Theorem \ref{Main}, we used Szemer\'{e}di's theorem. In fact, one can obtain a slightly better result by performing a more careful analysis. Let $N_1,N_2\geq 2$ be integers. Divide $[0,1]$ into $N_1$ equal pieces and then divide each of these small pieces into $N_2$  pieces of equal length. In total, we have $N_1N_2$ small pieces of length $(N_1N_2)^{-1}.$ These pieces can be grouped into $N_2$ many `arithmetic progressions' (more precisely speaking, intervals whose centres form an arithmetic progression) of length $N_1$ with gap $N_1^{-1}.$ To proceed further, we use the standard notation $r_k(N)$ to denote the largest cardinality of a subset of $\{1,\dots,N\}$ which contains no arithmetic progressions of length $k$.  We want to select a certain number of small pieces of length $(N_1N_2)^{-1}$ such that we do not get any  arithmetic progressions of length   of length $k$ with gap $N_1^{-1}.$ This number can be bounded above by $N_2r_k(N_1).$ For each $\epsilon>0,$ if $N_1$ is large enough then we can replace $r_k(N_1)$ by $\epsilon N_1.$ This is what we did in the proof of Theorem \ref{Main}. However, there are now stronger  quantitative upper bounds for $r_k(N)$ than the  Szemer\'{e}di bound of $\eps N$.  For example, we have the following result due to Gowers, see \cite{Go}.
\begin{thm}(Gowers)
	For each $k\geq 3,$ there are constants $c_k,C_k>0$ such that
	\[
	r_k(N)\leq C_k \frac{N}{(\log\log N)^{c_k}}
	\]
	for all $N\geq 100.$ Here $c_k$ can be chosen as $2^{-2^{k+9}}.$
\end{thm}
By applying this result, it should be possible to improve  Theorem \ref{Main} by replacing $\Delta^\alpha$ with $f(\Delta)$ for a suitable increasing function $f$. One can in principle compute $f$ precisely but we decided not to pursue the details. We suspect that $f(\Delta)=\Delta^{1/\log\log\log \Delta}$ will probably do the job. (This is obtained by replacing $\alpha$ in $(**)$ with $1/\log\log\log N$.)

Finally, we remark that the arguments in this paper provide a road map for translating estimates for $r_k(N)$ into statements of the type presented in Theorem \ref{Main}.  In fact, if one could establish sufficiently good estimates for $r_k(N)$, then one could prove the Erd\H{o}s  conjecture.

\subsection{Sharpness and non-integer sets}
Instead of trying to reduce the uncertainty, $f(\Delta)$, as discussed in the previous subsection, consider $f(\Delta) =  \Delta^{\alpha}$ for some  $\alpha\in (0,1)$, as in Theorem \ref{Main}.    What is interesting now is improving the `largeness' condition
\[
\#A\cap [0,n]\geq \frac{n}{(\log n)^{\gamma}}
\]
in the statement of Theorem \ref{Main}. In general, one can try to prove Theorem \ref{Main} with $n/(\log n)^{\gamma}$ being replaced by a general increasing function $g(n).$ A natural question to ask is whether $g(n)$ can be chosen to be $n^{\delta}$ for a $\delta>0.$ As we have seen in Section \ref{Mazur}, this is probably not true for $\delta=1/2.$ Also we suspect  that $\delta$ can depend on $\alpha.$

Our argument works not only for integer sets. Indeed, if we consider $A\subset\mathbb{R}^{+}$ and require that $A$ is uniformly $\delta$-discrete for a number $\delta>0$, i.e. $\inf_{a,b\in A} |a-b|>\delta>0,$ then all the arguments in the proof of Theorem \ref{Main} apply in this case and we have the same result. Clearly, in this non-integer case, one cannot hope to find exact arithmetic progressions since  small perturbations can destroy all exact progressions. Thus, the notion of `almost arithmetic progressions' we introduce here is very natural. In this case, one can try to lowering the `uncertainty' $f(\Delta)$ which was discussed in the previous paragraph.

	\section{Acknowledgements}
	JMF acknowledges  financial support  from  an EPSRC Standard Grant (EP/R015104/1).  HY was financially supported by the University of St Andrews. The authors thank Sam Chow for suggesting  Mazur's near miss problem.
	
	\providecommand{\bysame}{\leavevmode\hbox to3em{\hrulefill}\thinspace}
	\providecommand{\MR}{\relax\ifhmode\unskip\space\fi MR }
	\providecommand{\MRhref}[2]{%
		\href{http://www.ams.org/mathscinet-getitem?mr=#1}{#2}
	}
	\providecommand{\href}[2]{#2}


\begin{thebibliography}{dABCFS11}
	

\bibitem[Fr19]{fraserAP} J.~M. Fraser. \emph{Almost arithmetic progressions in the primes and other large sets}. Amer. Math. Monthly, (to appear).

%\bibitem[FSY19]{kota} J. M. Fraser, K. Saito and H. Yu.  \emph{Dimensions of sets which uniformly avoid arithmetic progressions}, Int. Math. Res. Not., available at  https://arxiv.org/abs/1705.03335.


\bibitem[FY18]{fraseryu} J. M. Fraser and H. Yu. \emph{Arithmetic patches, weak tangents, and dimension}, Bull. Lond. Math. Soc.,  \textbf{ 50}, 85--95, 2018.

\bibitem[FK78]{FK78} H. Furstenberg and Y. Katznelson. \emph{An ergodic Szemer\'edi theorem for commuting transformations}, J. Analyse Math., \textbf{34}, 275--291, 1978. 

\bibitem[G01]{Go}
T. Gowers, \emph{A new proof of Szemeredi's theorem}, Geom. Funct. Anal., \textbf{11}, 465--588, 2001.	

\bibitem[GT08]{greentao} B. Green and T. Tao. \emph{The primes contain arbitrarily long arithmetic progressions}. Ann. of Math.,   \textbf{167}, 481--547, 2008.
		
	
	
		
		
		\bibitem[M04]{MA} B. Mazur. \emph{Perturbations, Deformations, and Variations (and``Near-Misses") in Geometry, Physics, and Number Theory}. Bull. Amer. Math. Soc. (N.S.), \textbf{41}(3), 307-336, 2004.
		
		\bibitem[{S}75]{SZ}
		E. {Szemer\'edi}. \emph{{On sets of integers containing no $k$ elements in
				arithmetic progression.}}, {Acta Arith.}, \textbf{27}, 199--245, 1975.
			
		\bibitem[TZ15]{TZ15} T. Tao and T. Ziegler. \emph{A multi-dimensional Szemer\'edi theorem for the primes via a correspondence principle}. Israel J. Math., \textbf{207}, 203--228, 2015.
			
			\bibitem[vdW27]{VW}
		van~der Waerden. \emph{Beweis einer baudetschen Vermutung}, Nieuw Arch. Wisk, 212--216, 1927.
		
	
	\end{thebibliography}
\end{document}